\documentclass[12pt]{amsart}
\usepackage[cp1251]{inputenc}
\usepackage[english]{babel}
\usepackage{amsmath}
\usepackage{amssymb}
\usepackage{amsfonts}

\begin{document}

\title[On a general discrete boundary value problem]
      {On a General Discrete Boundary Value Problem }
\author[V. Vasilyev]{Vladimir Vasilyev}
\address{Chair of Applied Mathematics and Computer Modeling\\
    Belgorod State National Research University\\
         Pobedy street 85, Belgorod 308015, Russia}

        \email{vbv57@inbox.ru}

\author[A. Vasilyev]{Alexander Vasilyev}
\address{Chair of Applied Mathematics and Computer Modeling\\
    Belgorod State National Research University\\
         Pobedy street 85, Belgorod 308015, Russia}

        \email{alexvassel@gmail.com}

  \author[A. Mashinets]{Anastasia Mashinets}
\address{Chair of Applied Mathematics and Computer Modeling\\
    Belgorod State National Research University\\
         Pobedy street 85, Belgorod 308015, Russia}

        \email{anastasia.kho@yandex.ru}

 %\urladdr{http://math....}
%\thanks{Received , accepted .}
%\thanks{Research supported  under grant number...}
%\date{February 26, 2017}

\begin{abstract}
We study a general discrete boundary value problem in Sobolev--Slobodetskii spaces in a plane quadrant and reduce it to a system of integral equations. We show a solvability of the system for a small size of discreteness starting from a solvability of its continuous analogue.
\end{abstract}

\maketitle
%\tableofcontents

\section{Introduction}

The theory of pseudo-differential operators \cite{Ta,Tr} and related equations and boundary value problems \cite{E} exists more than a half-century, and up to now it takes attention of scholars. Also there are some discrete theories of boundary value problems for partial differential equations \cite{S,R}, but these studies are not applicable for pseudo-differential equations. According to this fact the first author has initiated a studying discrete theory of pseudo-differential equations \cite{VV1,V6,V3,V4} having in mind forthcoming studies their approximation properties and applications to computational algorithms \cite{VVT,MVV}.

Since model equations in \cite{E} were studied in a half-space, it is a canonical domain for manifold with smooth boundary, next step was done with a cone, it is a canonical domain for manifold with conical points at boundary \cite{V1}. At this step one needs a special factorization for an elliptic symbol, it was transferred to the discrete case \cite{V5}using analogies with Fourier series \cite{Ed}.

This paper is devoted to a model discrete pseudo-differential equation and discrete boundary value problem in a quadrant in the plane and its solvability starting from solvability their continuous analogues under small parameter of discreteness.

\section{Digital pseudo-differential operators and discrete equations}

Here we introduce some starting concepts and results which will help us moving to statement of a general boundary value problem.

Let $\mathbb Z^2$ be the integer lattice in a plane, $K=\{x\in\mathbb R^2: x=(x_1,x_2), x_1>0, x_2>0\}$ be the first quadrant, $K_d=h\mathbb Z^2\cap K, h>0$, $\mathbb T^2=[-\pi,\pi]^2, \hbar=h^{-1}$. We consider functions of a discrete variable $u_d(\tilde x), \tilde x=(\tilde x_1,\tilde x_2)\in h\mathbb Z^2$.

We use also notations $\zeta^2=h^{-2}((e^{ih\cdot\xi_1}-1)^2+(e^{ih\cdot\xi_2}-1)^2)$ and $S(h\mathbb Z^2)$ for the discrete analogue of the Schwartz space of infinitely differentiable rapidly decreasing at infinity functions.

{\bf Definition 1.} {\it
The space $H^s(h\mathbb Z^2)$ consists of discrete functions and it is a closure of the space  $S(h\mathbb Z^2)$ with respect to the norm
\begin{equation}\label{1}
||u_d||_s=\left(\int\limits_{\hbar\mathbb T^2}(1+|\zeta^2|)^s|\tilde u_d(\xi)|^2d\xi\right)^{1/2},
\end{equation}
where $\tilde u_d(\xi)$ denotes the discrete Fourier transform
$$
(F_du_d)(\xi)\equiv\tilde u_d(\xi)=\sum\limits_{\tilde x\in h\mathbb Z^2}e^{i\tilde x\cdot\xi}u_d(\tilde x)h^2,~~~\xi\in\hbar\mathbb T^2.
$$
}

Let $A_d(\xi)$ be a measurable periodic function defined in  $\mathbb R^2$ with the basic cube of periods $\hbar\mathbb T^2$.

{\bf Definition 2.} {\it
A digital pseudo-differential operator $A_d$ with the symbol $A_d(\xi)$ in discrete quadrant  $K_d$ is called the following operator
\[
(A_du_d)(\tilde x)=\sum\limits_{\tilde y\in h\mathbb Z^2}h^2\int\limits_{\hbar\mathbb T^2} A_d(\xi)e^{i(\tilde y-\tilde x)\cdot\xi}\tilde u_d(\xi)d\xi,~~~\tilde x\in K_d,
\]
}

Here we will consider symbols satisfying the condition
\[
c_1(1+|\zeta^2|)^{\alpha/2}\leq|A_d(\xi)|\leq c_2(1+|\zeta^2|)^{\alpha/2}
\]
with positive constants   $c_1, c_2$ non-depending on $h$. This class of symbols satisfying this condition will be denoted by $E_{\alpha}$. The number
 $\alpha\in\mathbb R$ is called an order of the digital pseudo-differential operator $A_d$.

We study solvability of the discrete equation
\begin{equation}\label{2}
(A_du_d)(\tilde x)=0,~~~\tilde x\in K_d,
 \end{equation}
in the space $H^s(K_d)$, and for this purpose we need certain specific domains of two-dimensional complex space  $\mathbb C^2.$ A domain of the type ${\mathcal T}_h(K)=\hbar\mathbb T^2+iK$ is called a tube domain over the quadrant $K$. We will work with holomorphic functions   $f(x+i\tau)$ in such domains ${\mathcal T}_h(K)$.

{\bf Definition 3.} {\it
Periodic wave factorization of the symbol
$A_d(\xi)\in E_{\alpha}$ is called its representation in the form
\[
A_d(\xi)=A_{d,\neq}(\xi)A_{d,=}(\xi),
\]
where the factors  $A_{d,\neq}(\xi), A_{d,=}(\xi)$ admit holomorphic continuation into tube domains ${\mathcal T}_h(K), {\mathcal T}_h(-K)$ respectively satisfying the estimates
\[
c_1(1+|\hat\zeta^2|)^{\frac{\ae}{2}}\leq|A_{d,\neq}(\xi+i\tau)|\leq c'_1(1+|\hat\zeta^2|)^{\frac{\ae}{2}},
\]
\[
c_2(1+|\hat\zeta^2|)^{\frac{{\alpha-\ae}}{2}}\leq|A_{d,=}(\xi-i\tau)|\leq c'_2(1+|\hat\zeta^2|)^{\frac{{\alpha-\ae}}{2}},
\]
with positive constants $c_1, c'_1, c_2, c'_2$ non-depending on $h$;
\[
\hat\zeta^2\equiv\hbar^2\left((e^{ih(\xi_1+i\tau_1)}-1)^2+(e^{ih(\xi_2+i\tau_2)}-1)^2\right),~~~\xi=(\xi_1,\xi_2)\in\hbar\mathbb T^2,
\]
\[
\tau=(\tau_1,\tau_2)\in K.
\]

The number $\ae\in{\mathbb R}$ is called an index of periodic wave factorization.
}

{\it Everywhere below we assume that we have this periodic wave factorization of the symbol $A_d(\xi)$ with the index $\ae$.}

Using methods developed in \cite{VV1} we can prove the following result.

{\bf Theorem 1.} {\it
Let $\ae-s=n+\delta, n\in\mathbb N, |\delta|<1/2$. Then a general solution of the equation  \eqref{2} has the following form
\begin{equation}\label{3}
\tilde u_d(\xi)= A^{-1}_{d,\neq}(\xi)\left(\sum\limits_{k=0}^{n-1}\left(\tilde c_k(\xi_1)\zeta_2^k+\tilde d_k(\xi_2)\zeta_1^k\right)\right),
\end{equation}
where  $\tilde c_k(\xi_1), \tilde d_k(\xi_2), k=0,1,\cdots,n-1,$ are arbitrary functions from-- $\widetilde H^{s_k}(h\mathbb T), s_k=s-\ae+k-1/2.$

The a priori estimate
\[
||u_d||_s\leq const\sum\limits_{k=0}^{n-1}([c_k]_{s_k}+[d_k]_{s_k}),
\]
holds, where $[\cdot]_{s_k}$ denotes a norm in the space    $H^{s_k}(h\mathbb T)$, and
$const$ doesn't depend on $h$.
}

\section{Discrete boundary value problem}

\subsection{Statement and solvability}

Starting from Theorem 1 we introduce the following boundary conditions
\begin{equation}\label{4}
\begin{aligned}
&(B_{d,j}u_d)(\tilde x_1,0)=b_{d,j}(\tilde x_1),
\\[3pt]
&(G_{d,j}u_d)(0,\tilde x_2)=g_{d,j}(\tilde x_2),~~~j=0,1.\cdots,n-1,
\end{aligned}
\end{equation}
where $B_{d,j}, G_{d,j}$ are digital pseudo-differential operators of order $\beta_j,\gamma_j\in{\mathbb R}$ with symbols $\widetilde B_{d,j}(\xi)\in E_{\beta_j}, \widetilde G_{d,j}(\xi)\in E_{\gamma_j}$
\[
\begin{aligned}
&(B_{d,j}u_d)(\tilde x)=\frac{1}{(2\pi)^2}\int\limits_{\hbar{\bf T}^2}\sum\limits_{\tilde y\in h{\bf Z}^2}e^{i\xi\cdot(\tilde x-\tilde y)}\widetilde B_{d,j}(\xi)\tilde u_d(\xi)d\xi,
\\[3pt]
&(G_{d,j}u_d)(\tilde x)=\frac{1}{(2\pi)^2}\int\limits_{\hbar{\bf T}^2}\sum\limits_{\tilde y\in h{\bf Z}^2}e^{i\xi\cdot(\tilde x-\tilde y)}\widetilde G_{d,j}(\xi)\tilde u_d(\xi)d\xi.
\end{aligned}
\]

One can rewrite boundary conditions
\eqref{4}
 in Fourier images
\begin{equation}\label{5}
\begin{aligned}
&\int\limits_{-\hbar\pi}^{\hbar\pi}\widetilde B_{d,j}(\xi_1,\xi_2)\tilde u_d(\xi_1,\xi_2)d\xi_2=\tilde b_{d,j}(\xi_1),
\\[3pt]
&\int\limits_{-\hbar\pi}^{\hbar\pi}\widetilde G_{d,j}(\xi_1,\xi_2)\tilde u_d(\xi_1,\xi_2)d\xi_1=\tilde g_{d,j}(\xi_2),~~~j=0,1.\cdots,n-1,
\end{aligned}
\end{equation}
so that according to properties of digital pseudo-differential operators and trace properties we need to require $b_{d,j}(\tilde x_1)\in H^{s-\beta_j-1/2}(h{\mathbb Z}), g_{d,j}(\tilde x_2)\in H^{s-\gamma_j-1/2}(h{\mathbb Z})$.

Multiplying the equality \eqref{3} by $\widetilde B_{d,j}(\xi_1,\xi_2)$ and $\widetilde G_{d,j}(\xi_1,\xi_2)$, integrating over $[-\hbar\pi,\hbar\pi]$ on $\xi_2$ and $\xi_1$, taking into account the conditions \eqref{5} we obtain the following $(2n\times 2n)$-system of linear integral equations
\begin{equation}\label{6}
\begin{aligned}
&\sum\limits_{k=0}^{n-1}\left(r_{jk}(\xi_1)\tilde c_k(\xi_1)+\int\limits_{-\hbar\pi}^{\hbar\pi}l_{jk}(\xi_1,\xi_2)\tilde d_k(\xi_2)d\xi_2\right)=\tilde b_{d,j}(\xi_1)
\\[3pt]
&\sum\limits_{k=0}^{n-1}\left(\int\limits_{-\hbar\pi}^{\hbar\pi}m_{jk}(\xi_1,\xi_2)\tilde c_k(\xi_1)d\xi_1+p_{jk}(\xi_2)\tilde d_k(\xi_2)\right)=\tilde g_{d,j}(\xi_2),
\\[3pt]
&j=0,1,\dots,n-1,
\end{aligned}
\end{equation}
with unknowns functions $\tilde c_k, \tilde d_k, k=0,1,\dots,n-1$. We have used the following notations
\[
r_{jk}(\xi_1)=\int\limits_{-\hbar\pi}^{\hbar\pi}\tilde B_{d,j}(\xi)A_{d,\neq}^{-1}(\xi)\zeta_2^kd\xi_2,~~~p_{jk}(\xi_2)=\int\limits_{-\hbar\pi}^{\hbar\pi}\tilde G_{d,j}(\xi)A_{d,\neq}^{-1}(\xi)\zeta_1^kd\xi_1,
\]
\[
l_{jk}(\xi_1,\xi_2)=\tilde B_{d,j}(\xi)A_{d,\neq}^{-1}(\xi)\zeta_1^k,~~~m_{jk}(\xi_1,\xi_2)=\tilde G_{d,j}(\xi)A_{d,\neq}^{-1}(\xi)\zeta_2^k,
\]
$j,k=0,1,\dots,n-1$.

Thus, we can formulate the following assertion.

{\bf Theorem 2.} {\it
The boundary value problem \eqref{2},\eqref{4} is uniquely solvable in the space $H^s(K_d)$ with data $b_{d,j}\in H^{s-\beta_j-1/2}(h\mathbb Z_+, g_{d,j}\in H^{s-\gamma_j-1/2}(h\mathbb Z_+$ if and only if the system \eqref{6} has the unique solution $\tilde c_k,\tilde d_k\in\tilde H^{s_k}(\hbar\mathbb T)$, $j,k=0,1,\dots,n-1$.
}

\subsection{Continuous case}

Here we will describe continuous boundary value problem which is related to considered discrete boundary value problem \eqref{2},\eqref{4}.

Let $A$ be a pseudo-differential operator
\[
(Au)(x)=\int\limits_{\mathbb R^2}\int\limits_{\mathbb R^2}\tilde A(\xi)e^{i\xi(y-x)}u(y)dyd\xi
\]
with symbol $\tilde A(\xi)$ satisfying the condition
\begin{equation}\label{7}
|\tilde A(\xi)|\sim(1+|\xi|)^{\alpha}
\end{equation}
and admitting the wave factorization with respect to $K$
\[
\tilde A(\xi)=A_{\neq}(\xi)\cdot A_=(\xi).
\]
with index $\ae$ such that $\ae-s=n+\delta, n\in\mathbb N, |\delta|<1/2$.

Further, let $B_j,G_j, j=0,1,\dots,n-1$ be pseudo-differential operators with symbols $\tilde B_j(\xi),\tilde G_j(\xi)$ satisfying the condition \eqref{7} with $\beta_j,\gamma_j$ instead of $\alpha$.

The following boundary value problem
\begin{equation}\label{8}
\begin{aligned}
&(Au)(x)=0,~~~x\in K,
\\[3pt]
&(B_ju)(x_1,0)=b_j(x_1),
\\[3pt]
&(G_ju)(0,x_2)=g_j(x_2),~~~j=0,1,\dots,n-1
\end{aligned}
\end{equation}
is a continuous analogue of the discrete boundary value problem \eqref{2},\eqref{4}. It was shown in \cite{V1} the problem \eqref{8} is equivalent to the following system of integral equations
\begin{equation}\label{9}
\begin{aligned}
&\sum\limits_{k=0}^{n-1}\left(R_{jk}(\xi_1)\tilde C_k(\xi_1)+\int\limits_{-\infty}^{+\infty}L_{jk}(\xi_1,\xi_2)\tilde D_k(\xi_2)d\xi_2\right)=\tilde b_{j}(\xi_1)
\\[3pt]
&\sum\limits_{k=0}^{n-1}\left(\int\limits_{-\infty}^{+\infty}M_{jk}(\xi_1,\xi_2)\tilde C_k(\xi_1)d\xi_1+P_{jk}(\xi_2)\tilde D_k(\xi_2)\right)=\tilde g_{j}(\xi_2),
\\[3pt]
&j=0,1,\dots,n-1,
\end{aligned}
\end{equation}
with unknowns functions $\tilde C_k, \tilde D_k, k=0,1,\dots,n-1$. The following notations
are used
\[
R_{jk}(\xi_1)=\int\limits_{-\infty}^{+\infty}\tilde B_{j}(\xi)A_{\neq}^{-1}(\xi)(i\xi_2)^kd\xi_2,~~~P_{jk}(\xi_2)=\int\limits_{-\infty}^{+\infty}\tilde G_{j}(\xi)A_{\neq}^{-1}(\xi)(i\xi_1)^kd\xi_1,
\]
\[
L_{jk}(\xi_1,\xi_2)=\tilde B_{j}(\xi)A_{\neq}^{-1}(\xi)(i\xi_1)^k,~~~M_{jk}(\xi_1,\xi_2)=\tilde G_{j}(\xi)A_{\neq}^{-1}(\xi)(i\xi_2)^k,
\]
$j,k=0,1,\dots,n-1$. If we can solve the system \eqref{9} and find $\tilde C_k, \tilde D_k, k=0,1,\dots,n-1$ the the solution of the boundary value problem \eqref{9} can be constructed by the formula \cite{V1}
\begin{equation}\label{10}
\tilde u(\xi)= A^{-1}_{\neq}(\xi)\left(\sum\limits_{k=0}^{n-1}\left(\tilde C_k(\xi_1)(i\xi_2)^k+\tilde D_k(\xi_2)(i\xi_1)^k\right)\right),
\end{equation}
where  $\tilde C_k(\xi_1), \tilde D_k(\xi_2), k=0,1,\cdots,n-1,$ are arbitrary functions from $\widetilde H^{s_k}(\mathbb R), s_k=s-\ae+k-1/2$.

Our next problems are the following. Given operator $A$ and boundary operators $B_j,G_j$ how to choice the digital operators $A_d$ and $B_{d,j},G_{d,j}$ to obtain the implication: the unique solvability of the system \eqref{9} gives the unique solvability of the system \eqref{6} for enough small $h$. This question will be discussed in the next section.

\section{Comparison theorems}

\subsection{Projection method}

Let us introduce the following space of vector-functions:
\[
\tilde{\bf H}^{\Lambda}(\mathbb R)=\tilde{\bf H}^S(\mathbb R)\oplus\tilde{\bf H}^S(\mathbb R),~~~\tilde{\bf H}^S(\mathbb R)=\oplus\sum\limits_{k=0}^{n-1}\tilde{H}^{s_k}(\mathbb R),
\]

Norms in these spaces will be defined in the following way. For $f\in\tilde{\bf H}^{S}(\mathbb R), f=(f_0,\dots,f_{n-1}), f_k\in\tilde{H}^{s_k}(\mathbb R), g\in\tilde{\bf H}^{S}(\mathbb R), g=(g_0,\dots,g_{n-1}), g_k\in\tilde{H}^{s_k}(\mathbb R)$ we put
\[
||f||_{S}=\sum\limits_{k=0}^{n-1}||f_k||_{s_k}.~~~||g||_{S}=\sum\limits_{k=0}^{n-1}||g_k||_{s_k},
\]
and if $F\in\tilde{\bf H}^{\Lambda}(\mathbb R), F=(f,g),  f\in\tilde{\bf H}^{S}(\mathbb R), g\in\tilde{\bf H}^{S}(\mathbb R)$ we put
\[
||F||_{\Lambda}=||f||_{S}+||g||_{S}.
\]
%A norm in the space $\tilde{\bf H}^{\Lambda}(\mathbb R)$ is defined in the same way.

Let us introduce the following notations. We denote the system \eqref{9} in the following way
\[
\begin{pmatrix}
R&L\\
M&P
\end{pmatrix}
\begin{pmatrix}
C\\
D
\end{pmatrix}=\begin{pmatrix}
B\\
G
\end{pmatrix},
\]
where $C=(\tilde c_0,\dots,\tilde c_{n-1})^T, D=(\tilde d_0,\dots,\tilde d_{n-1})^T, B=(\tilde b_0,\dots,\tilde b_{n-1})^T, G=(\tilde g_0,\dots,\tilde g_{n-1})^T$; operators $R,L,M,P$ acting in the space $\tilde{\bf H}^S(\mathbb R)$ are the following: $R$ is multiplier by the matrix-function $(r_{jk})_{j,k=0}^{n-1}$, $P$ is multiplier by the matrix-function $(p_{jk})_{j,k=0}^{n-1}$, $L$, $M$ are matrix integral operators with kernels $L_{jk}, M_{jk}$ respectively.

Further, we will denote $\Xi_h$ the restriction operator on the segment $\hbar\mathbb T$ so that  for $f\in\tilde{\bf H}^S(\mathbb R), f=(f_0,\dots,f_{n-1})$ the notation $\Xi_h f$ means the following
\[
\Xi_h f=(\chi_h f_0,\dots,\chi_h f_{n-1}),
\]
where $\chi_h$ is an indicator of $\hbar\mathbb T$.

We denote by $Q$ the operator
\[
Q=\begin{pmatrix}
R&L\\
M&P
\end{pmatrix}
\]

{\bf Theorem 3.} {\it Let $s-\beta_j>1, s-\gamma_j>2, j=0,1,\dots,n-1$.
We have the estimate
\[
||\Xi_hQ-Q\Xi_h||_{\tilde{\bf H}^{\Lambda}(\mathbb R)\rightarrow\tilde{\bf H}^{\Lambda}(\mathbb R)}\leq~const~ h^{\varepsilon},
\]
where
$$
\varepsilon =\min\limits_{0\leq j\leq n-1}\{s-\beta_j-1,s-\gamma_j-1\},
$$
const does not depend on $h, s_k=s-\ae+k-1/2, k=0,1,\dots,n-1$.
}

\begin{proof}
Obviously, the matrices $R,P$ give vanishing result in the norm, and we need to work with integral operators only. Let us consider the operator $L$, and extract one its component $L_{jk}$,
\[
\int\limits_{-\infty}^{+\infty}L_{jk}(\xi_1,\xi_2)\tilde D_k(\xi_2)d\xi_2,~~~L_{jk}(\xi_1,\xi_2)=\tilde B_{j}(\xi)A_{\neq}^{-1}(\xi)\xi_1^k.
\]

We have
\[
\chi_h(\xi_1)\int\limits_{-\infty}^{+\infty}L_{jk}(\xi_1,\xi_2)\tilde D_k(\xi_2)d\xi_2-
\int\limits_{-\bar\pi}^{+\hbar\pi}L_{jk}(\xi_1,\xi_2)\tilde D_k(\xi_2)d\xi_2
\]
\[
=\begin{cases}
\left(\int\limits_{-\infty}^{-\hbar\pi}+\int\limits_{\hbar\pi}^{+\infty}\right)L_{jk}(\xi_1,\xi_2)\tilde D_k(\xi_2)d\xi_2,~~~\xi_1\in\hbar\mathbb T,\\
-\int\limits_{-\hbar\pi}^{\hbar\pi}L_{jk}(\xi_1,\xi_2)\tilde D_k(\xi_2)d\xi_2,~~~\xi_1\notin\hbar\mathbb T.
\end{cases}
\]

Let us consider the first case and estimate
\[
\left|\int\limits_{\hbar\pi}^{+\infty}L_{jk}(\xi_1,\xi_2)\tilde D_k(\xi_2)d\xi_2\right|\leq~const\int\limits_{\hbar\pi}^{+\infty}(1+|\xi|)^{\beta_j-\ae}|\xi_1|^k|\tilde D_k(\xi_2)|d\xi_2
\]
\[
\leq~const\int\limits_{\hbar\pi}^{+\infty}(1+|\xi|)^{\beta_j-s+1/2}|\tilde D_k(\xi_2)|(1+|\xi_2|)^{s_k}d\xi_2
\]
(we have taken into account $s_k=s-\ae+k-1/2$ and now we apply the Cauchy--Schwartz inequality)
\[
\leq~const(1+|\xi_1|+\hbar)^{\beta_j-s+1}||\tilde D_k||_{s_k}\leq~const~h^{s-\beta_j-1}||D_k||_{s_k}.
\]

Squaring the latter inequality, multiplying by $(1+|\xi_|)^{2s_k}$ and integrating over $\hbar\mathbb T$ we obtain
\[
\int\limits_{-\hbar\pi}^{\hbar\pi}(1+|\xi_|)^{2s_k}\left|\int\limits_{\hbar\pi}^{+\infty}L_{jk}(\xi_1,\xi_2)\tilde D_k(\xi_2)d\xi_2\right|^2d\xi_1
\]
\[
\leq~const~h^{2(s-\beta_j-1)}||D_k||^2_{s_k}\int\limits_{0}^{+\infty}(1+|\xi_|)^{2s_k}d\xi_1\leq~const~h^{2(s-\beta_j-1)}||D_k||^2_{s_k}
\]

For the second case ($|\xi_1|>\hbar\pi$) we obtain
\[
\left|\int\limits_{-\hbar\pi}^{+\hbar\pi}L_{jk}(\xi_1,\xi_2)\tilde D_k(\xi_2)d\xi_2\right|\leq~const\int\limits_{-\hbar\pi}^{+\hbar\pi}(1+|\xi|)^{\beta_j-\ae}|\xi_1|^k|\tilde D_k(\xi_2)|d\xi_2
\]
\[
\leq~const\int\limits_{-\hbar\pi}^{+\hbar\pi}(1+|\xi|)^{\beta_j-\ae}|\xi_1|^{n-1}(1+|\xi_2|)^{-s_k}|\tilde D_k(\xi_2)|(1+|\xi_2|)^{s_k}d\xi_2
\]
\[
\leq~const~|\xi_1|^{n-1}(1+|\xi_1|)^{-s_k}\int\limits_{-\hbar\pi}^{+\hbar\pi}(1+|\xi|)^{\beta_j-\ae}|\tilde D_k(\xi_2)|(1+|\xi_2|)^{s_k}d\xi_2
\]
(we apply the Cauchy--Schwartz inequality in the integral)
\[
\leq~const~(1+|\xi_1|)^{n-s_k-1}(1+|\xi_1|)^{\beta_j-\ae+1/2k}||D_k||_{s_k}
\]

Squaring the latter inequality, multiplying by $(1+|\xi_|)^{2s_k}$ and integrating over $\mathbb R\setminus\hbar\mathbb T$ we obtain
\[
\left(\int\limits_{-\infty}^{-\hbar\pi}+\int\limits_{\hbar\pi}^{+\infty}\right)(1+|\xi_|)^{2s_k}\left|\int\limits_{\hbar\pi}^{+\infty}L_{jk}(\xi_1,\xi_2)\tilde D_k(\xi_2)d\xi_2\right|^2d\xi_1
\]
\[
\leq~const||D_k||^2_{s_k}\int\limits_{\hbar\pi}^{+\infty}(1+\xi_1)^{2n-2+2\beta_j+1-2\ae}d\xi_1\leq~const||D_k||^2_{s_k}h^{2s-2\beta_j+2\delta},
\]
since $2n+2\beta_j-2\ae=2n+2\beta_j-2(s+n+\delta)=2\beta_j-2s-2\delta<0$.

Thus, we have proved that
\[
||\chi_hL_{jk}-L_{jk}\chi_h||_{H^{s_k}(\mathbb R)\rightarrow H^{s_k}(\mathbb R)}\leq~const~h^{s-\beta_j-1},
\]
since $s-\beta_j-1<s-\beta_j+\delta$.

Almost the same inequality can be obtained for $M_{jk}$
\[
||\chi_hM_{jk}-M_{jk}\chi_h||_{H^{s_k}(\mathbb R)\rightarrow H^{s_k}(\mathbb R)}\leq~const~h^{s-\gamma_j-1},
\]
These estimates complete the proof.
\end{proof}

{\bf Corollary 1.} {\it Under conditions of Theorem 3 the invertibility of the operator $Q$ in the space $\tilde{\bf H}^{\Lambda}(\mathbb R)$ implies the invertibility of the operator $\Xi_hQ\Xi_h$ in the space $\tilde{\bf H}^{\Lambda}(\hbar\mathbb T)$ for enough small $h$.

}

\begin{proof}
We apply the results of the paper \cite{K} which imply the following. If
\[
||\Xi_hQ-Q\Xi_h||_{\tilde{\bf H}^{\Lambda}(\mathbb R)\rightarrow\tilde{\bf H}^{\Lambda}(\mathbb R)}\rightarrow 0,~~~h\to 0
\]
then the equation in the space $\tilde{\bf H}^{\Lambda}(\mathbb R)$
\begin{equation}\label{11}
Qu=v
\end{equation}
admits applying so called {\it projection method}. In other words it means that unique solvability of the equation \eqref{11} in the space $\tilde{\bf H}^{\Lambda}(\mathbb R)$ implies unique solvability of the equation
\begin{equation}\label{12}
\Xi_hQ\Xi_hu=\Xi_hv
\end{equation}
in the space $\tilde{\bf H}^{\Lambda}(\hbar\mathbb T)$ for enough small $h$. Moreover, if there is bounded operator $Q^{-1}$ in the space $\tilde{\bf H}^{\Lambda}(\mathbb R)$ then there is bounded operator $(\Xi_hQ\Xi_h)^{-1}$ for enough small $h$ and
\[
||(\Xi_hQ\Xi_h)^{-1}||_{\tilde{\bf H}^{\Lambda}(\hbar\mathbb T)\rightarrow\tilde{\bf H}^{\Lambda}(\hbar\mathbb T)}\leq~const,
\]
where const does not depend on $h$.

Indeed, a reader can easily verify that
\[
||(\Xi_hQ\Xi_h)^{-1}-\Xi_hQ^{-1}\Xi_h||_{\tilde{\bf H}^{\Lambda}(\hbar\mathbb T)\rightarrow\tilde{\bf H}^{\Lambda}(\hbar\mathbb T)}\rightarrow 0,~~~h\to 0.
\]
\end{proof}

\subsection{Discrete and continuous}

To compare discrete and continuous operators we need a special choice of discrete operators. We will do it in the following way.

The symbol $A_d(\xi)$ of the discrete operator  $A_d$ will be constructed as follows. Given wave factorization for  $\tilde A(\xi)$
\[
\tilde A(\xi)=A_{\neq}(\xi)\cdot A_=(\xi)
\]
 we take restrictions of factors  $A_{\neq}(\xi), A_=(\xi)$ on $\hbar\mathbb T^2$ and periodically continue them into  $\mathbb R^2$. We denote these elements by $A_{d,\neq}(\xi), A_{d,=}(\xi)$ and construct the periodic symbol  $A_d(\xi)$ which admits periodic wave factorization with respect to  $K$
\[
A_d(\xi)=A_{d,\neq}(\xi)\cdot A_{d,=}(\xi)
\]
with the same index $\ae$. We construct discrete pseudo-differential operators $B_{d,j}, G_{d,j}$ taking their symbol as restrictions of symbols $\tilde B_j(\xi), \tilde G_j(\xi)$ on $\hbar\mathbb T^2$ with periodical continuations into $\mathbb R^2,, j=0,1,\dots,n-1$. The discrete boundary functions $b_{d,j}$, $g_{d,j}$ are constructed in the same way.
Thus, we have the corresponding discrete boundary value problem \eqref{2},\eqref{4}.

{\bf Lemma 1.} {\it The estimate
\[
|(i\xi_m)^k-\zeta_m^k|\leq~const~h|\xi_m|^{k+1}
\]
holds for $\xi_m\in\hbar\mathbb T, m=1,2$ const does not depend on $h$.
}

\begin{proof}

First, we estimate
\[
|\zeta_1|=\left|\sum\limits_{\nu=1}^{\infty}\frac{(i\xi_1)^{\nu+1}h^{\nu}}{(\nu+1)!}\right|=|\xi_1|\left|\sum\limits_{\nu=0}^{\infty}\frac{(i\xi_1)^{\nu}h^{m}}{(\nu+1)!}\right|
\leq|\xi_1|\sum\limits_{\nu=0}^{\infty}\frac{(|\xi_1|h)^{\nu}}{\nu!}
\]
\[
=|\xi_1|e^{|\xi_1|h}\leq|\xi_1|e^{\pi}.
\]

Second,
\[
|\zeta_1-i\xi_1|=|\hbar(e^{i\xi_1h}-1)-i\xi_1|=\left|\sum\limits_{\nu=1}^{\infty}\frac{(i\xi_1)^{\nu+1}h^{\nu}}{(\nu+1)!}\right|
\]
\[
\leq|\xi_1|^2h\sum\limits_{\nu=0}^{\infty}\frac{|\xi_1|^{\nu}h^{\nu}}{\nu!}=|\xi_1|^2he^{|\xi_1|h}\leq|\xi_1|^2he^{\pi}.
\]

We have
\[
\zeta_1^k-(i\xi_1)^k=(\zeta_1-i\xi_1)\left(\sum\limits_{\nu=0}^{k-1}\zeta_1^\nu(i\xi_1)^{k-1-\nu}\right),
\]
and thus
\[
|\zeta_1^k-(i\xi_1)^k|\leq|\zeta_1-i\xi_1|\sum\limits_{\nu=0}^{k-1}|\zeta_1|^\nu|\xi_1|^{k-1-\nu}
\]
Applying above estimates we obtain required inequality.

\end{proof}

{\bf Lemma 2.} {\it Let $s-\beta_j>2, s-\gamma_j>2, j=0,1,\dots,n-1$.The following estimates
\[
\begin{aligned}
&|L_{jk}(\xi_1,\xi_2)-l_{jk}(\xi_1,\xi_2)|\leq~const~h(1+|\xi|)^{\beta_j-\ae+k+1},
\\[3pt]
&|M_{jk}(\xi_1,\xi_2)-m_{jk}(\xi_1,\xi_2)|\leq~const~h(1+|\xi|)^{\gamma_j-\ae+k+1},
\\[3pt]
&|R_{jk}(\xi_1)-r_{jk}(\xi_1)|\leq~const~h(1+|\xi_1|)^{\beta_j-\ae+k+2},
\\[3pt]
&|P_{jk}(\xi_2)-p_{jk}(\xi_2)|\leq~const~h(1+|\xi_1|)^{\gamma_j-\ae+k+2}
\end{aligned}
\]
hold for $\xi_1,\xi_2\in\hbar\mathbb T$.
}

\begin{proof}
According to above conventions we have for $\xi\in\hbar\mathbb T^2$ and using Lemma 1 we have
\[
|L_{jk}(\xi_1,\xi_2)-l_{jk}(\xi_1,\xi_2)|=|\tilde B_j(\xi)A_{\neq}^{-1}(\xi)-\tilde B_{d,j}\xi)A_{d,\neq}^{-1}(\xi)||B_j(\xi)||\xi_1^k-\zeta_1^k|
\]
\[
\leq~const~(1+|\xi|)^{\beta_j-\ae}h|\xi_1|^{k+1}\leq~const~h(1+|\xi|)^{\beta_j-\ae+k+1}.
\]

Further,
\[
|R_{jk}(\xi_1)-r_{jk}(\xi_1)|=\left|\int\limits_{-\infty}^{+\infty}\tilde B_j(\xi)A_{\neq}^{-1}(\xi)\xi_2^kd\xi_2-\int\limits_{-\hbar\pi}^{+\hbar\pi}\tilde B_{d,j}(\xi)A_{d,\neq}^{-1}(\xi)\zeta_2^kd\xi_2\right|
\]
\[
\leq\int\limits_{-\hbar\pi}^{+\hbar\pi}|\tilde B_{d,j}(\xi)A_{d,\neq}^{-1}(\xi)||\xi_2^k-\zeta_2^k|d\xi_2+\left(\int\limits_{-\infty}^{-\hbar\pi}+\int\limits_{\hbar\pi}^{+\infty}\right)
|\tilde B_j(\xi)A_{\neq}^{-1}(\xi)\xi_2^k|d\xi_2.
\]
For the first integral we have
\[
\int\limits_{-\hbar\pi}^{+\hbar\pi}|\tilde B_{d,j}(\xi)A_{d,\neq}^{-1}(\xi)||\xi_2^k-\zeta_2^k|d\xi_2\leq~const~h\int\limits_{-\hbar\pi}^{+\hbar\pi}(1+|\xi|)^{\beta_j-\ae+k+1}d\xi_2
\]
\[
\leq~const~h(1+|\xi_1|)^{\beta_j-\ae+k+2},
\]
since $\beta_j-\ae+k+2<0, s-\beta_j>2.$

The second summand
\[
\left|\left(\int\limits_{-\infty}^{-\hbar\pi}+\int\limits_{\hbar\pi}^{+\infty}\right)
|\tilde B_j(\xi)A_{\neq}^{-1}(\xi)\xi_2^k|d\xi_2\right|\leq~const\int\limits_{\hbar\pi}^{+\infty}(1+|\xi_1|+|\xi_2|)^{\beta_j-\ae+k}d\xi_2
\]
\[
\leq~const(1+|\xi_1|+\hbar)^{\beta_j-\ae+k+1}\leq~const~h(1+|\xi_1|)^{\beta_j-\ae+k+2}.
\]

The same estimates are valid for $M_{jk}-m_{jk}$ and $P_{jk}$ with $\gamma_j$ instead of $\beta_j$.
\end{proof}

We introduce similar notations for the system \eqref{6} so that this system takes the form
\[
\begin{pmatrix}
r&l\\
m&p
\end{pmatrix}
\begin{pmatrix}
c\\
d
\end{pmatrix}=\begin{pmatrix}
B_d\\
G_d
\end{pmatrix},
\]
where
\[
q=\begin{pmatrix}
r&l\\
m&p
\end{pmatrix}
\]
is linear bounded operator acting in the space $\tilde{\bf H}^{\Lambda}(\hbar\mathbb T)$.

{\bf Theorem 4.} {\it Let $s-\beta_j>3, s-\gamma_j>3, j=0,1,\dots,n-1$. A comparison between operators $Q$ and $q$ is given by the estimate
\[
||\Xi_hQ\Xi_h-q||_{\tilde{\bf H}^{\Lambda}(\hbar\mathbb T)\rightarrow\tilde{\bf H}^{\Lambda}(\hbar\mathbb T)}\leq~const!h,
\]
where const does not depend on $h$.
}

\begin{proof}
We need to estimate $H^{s_k}(\hbar\mathbb T$-norms of the following elements
\[
(R_{jk}(\xi_1)-r_{jk}(\xi_1))f(\xi_1),~~~(R_{jk}(\xi_2)-p_{jk}(\xi_2))f(\xi_2),
\]
\[
\int\limits_{-\hbar\pi}^{\hbar\pi}(L_{jk}(\xi_1,\xi_2)-l_{jk}(\xi_1,\xi_2))f(\xi_2)d\xi_2,
\]
\[
\int\limits_{-\hbar\pi}^{\hbar\pi}(M_{jk}(\xi_1,\xi_2)-m_{jk}(\xi_1,\xi_2))f(\xi_1d\xi_1.
\]

We have according to Lemma 2
\[
|(R_{jk}(\xi_1)-r_{jk}(\xi_1))f(\xi_1)|\leq~const~h(1+|\xi_1|)^{\beta_j-\ae+k+2}|f(\xi_1)|.
\]
%\[
%=const~h(1+|\xi_1|)^{\beta_j-s+5/2}|f(\xi_1)|(1+|\xi_1|)^{s_k}
%\]
%since $\beta_j-\ae+k+2-s_k=-\beta_j-s+5/2$.
Multiplying the latter inequality by $(1+|\xi_1|)^{s_k}$, squaring, integrating over $\hbar\mathbb T$ and applying the Cauchy--Schwartz inequality we obtain
\[
\int\limits_{-\hbar\pi}^{+\hbar\pi}(1+|\xi_1|)^{2s_k}|R_{jk}(\xi_1)-r_{jk}(\xi_1)|^2|f(\xi_1)|^2d\xi_1
\leq~const~h^2||f||^2_{s_k}
\]
since $\beta_j-\ae+k+2<0$.

Let us consider
\[
\int\limits_{-\hbar\pi}^{\hbar\pi}(L_{jk}(\xi_1,\xi_2)-l_{jk}(\xi_1,\xi_2))f(\xi_2)d\xi_2.
\]
Using Lemma 2 we have
\[
\left|\int\limits_{-\hbar\pi}^{\hbar\pi}(L_{jk}(\xi_1,\xi_2)-l_{jk}(\xi_1,\xi_2))f(\xi_2)d\xi_2\right|
\]
\[
\leq~const~h\int\limits_{-\hbar\pi}^{\hbar\pi}(1+|\xi|)^{\beta_j-\ae+k+1}|f(\xi_2)|d\xi_2
\]
\[
\leq~const~h\int\limits_{-\hbar\pi}^{\hbar\pi}(1+|\xi|)^{\beta_j-s+5/2}|f(\xi_2)|(1+|\xi_2|)^{s_k}d\xi_2
\]
since $\beta_j-\ae+k+2-s_k=-\beta_j-s+5/2$. Now applying Cauchy--Schwartz inequality we find
\[
\left|\int\limits_{-\hbar\pi}^{\hbar\pi}(L_{jk}(\xi_1,\xi_2)-l_{jk}(\xi_1,\xi_2))f(\xi_2)d\xi_2\right|
\]
\[
\leq~const~h||f||_{s_k}\left(\int\limits_{-\hbar\pi}^{\hbar\pi}(1+|\xi_1|+|\xi_2|)^{2\beta_j-2s+5}d\xi_2\right)^{1/2}
\]
\[
\leq~const~h||f||_{s_k}(1+|\xi_1|)^{\beta_j-s+3}
\]
according to the condition $s-\beta_j>3$. Squaring, multiplying by $(1+|\xi_1|)^{2s_k}$ and integrating over $\hbar\mathbb T$ we conclude
\[
\int\limits_{-\hbar\pi}^{\hbar\pi}(1+|\xi_1|)^{2s_k}\left|\int\limits_{-\hbar\pi}^{\hbar\pi}
(L_{jk}(\xi_1,\xi_2)-l_{jk}(\xi_1,\xi_2))f(\xi_2)d\xi_2\right|^2d\xi_1
\]
\[
\leq~const~h^2||f||^2_{s_k}\int\limits_{-\hbar\pi}^{\hbar\pi}(1+|\xi_1|)^{2(\beta_j-s+3+s_k)}d\xi_1\leq~const~h^2||f||^2_{s_k}
\]
since $2(\beta_j-s+3+s_k)<-1$. Indeed, $2(\beta_j-s+3+s_k)=2(\beta_j-s+3+s-\ae+k-1/2)=2(\beta_j-s-\delta).$ Obviously, the inequality $2(\beta_j-s-\delta)<-1$ is equivalent to $s-\beta_j>-1-\delta$.

\end{proof}

{\bf Corollary 2.} {\it Under conditions of Theorem 4 the invertibility of the operator $Q$ in the space $\tilde{\bf H}^{\Lambda}(\mathbb R)$ implies the invertibility of the operator $q$ in the space $\tilde{\bf H}^{\Lambda}(\hbar\mathbb T)$ for enough small $h$.

}

\begin{proof}
Indeed, we have  the invertibility of $\Xi_hQ\Xi_h$ by Corollary 1 and the invertibility of $q$ is obtained by Theorem 4.

\end{proof}

\section*{Conclusion}

Main goal of the paper was to prove unique solvability of discrete boundary value problem having in main unique solvability of its continuous analogue. It was done by a special choice od discrete operator and discrete boundary conditions. We hope that estimates of Theorem 3,4 will help us to obtain some estimates for discrete and continuous solutions.

%\vspace{5cm}

\end{document}